\documentclass[12pt, a4paper]{article}

\usepackage[margin=15mm]{geometry}
\usepackage[square,numbers]{natbib}
\usepackage[width=0.88\textwidth]{caption}
\usepackage{amsfonts,mathtools,booktabs,array,multirow,float,url}

\providecommand{\keywords}[1]
{
	\small	
	\textbf{Keywords: } #1
}

\title{Stability switching in Lotka\textendash Volterra and Ricker\textendash type predator\textendash prey systems with arbitrary step size}

\author{Shamika Kekulthotuwage Don$^{1}$, Kevin Burrage$^{1,2}$, Kate Hemlstedt$^{1,2}$,\\ Pamela Burrage$^{1,2}$\\
	\footnotesize \textit{$^{1}$ School of Mathematical Sciences, Queensland University of Technology,  Brisbane, Australia}  \\
	\footnotesize \textit{$^{2}$ Centre for Data Science, Queensland University of Technology, Brisbane, Australia}
}

\begin{document} 

\maketitle	

\abstract{Dynamical properties of numerically approximated discrete systems may become inconsistent with those of the corresponding continuous\textendash time system. We present a qualitative analysis of the dynamical properties of two species Lotka\textendash Volterra and Ricker\textendash type predator\textendash prey systems under discrete and continuous settings. By creating an arbitrary time discretisation, we obtain stability conditions that preserve the characteristics of continuous\textendash time models and their numerically approximated systems. Here, we show that even small changes to some of the model parameters may alter the system dynamics unless an appropriate time discretisation is chosen to return similar dynamical behaviour observed in the corresponding continuous\textendash time system. We also found similar dynamical properties of the Ricker\textendash type predator\textendash prey systems under certain conditions. Our results demonstrate the need for preliminary analysis to identify which dynamical properties of approximated discretised systems agree or disagree with the corresponding continuous\textendash time systems. }

		
\section{Introduction}

Population density variations of interacting species in ecosystems are modelled as discrete\textendash time systems. Difference equations, which consider time in discrete form (such as the Ricker model \citep{ricker1958handbook}) are often applied in ecological modelling to predict new values from those evaluated at a previous discrete time step. Our study is motivated by the dynamical properties of discrete systems with arbitrary step size. We begin with mentioning a simple Ricker\textendash type predator\textendash prey system in discrete form, with unit time step, given by
\begin{gather}\label{eq:Baxter}
	\begin{split}
		N(t+1)&=N(t)e^{r(1-\frac{N(t)}{K})-\alpha P(t)}\\
		P(t+1)&=P(t)e^{\alpha \gamma N(t)-c},
	\end{split}
\end{gather}
where $N(t)$ and $P(t)$ denote the prey and predator population at time $t$ respectively, $r$ is the rate of prey population increase, $K$ is prey carrying capacity, $\alpha$ is the predator attack rate, $\gamma$ is the conversion rate of eaten prey to sustenance for the predators and $c$ is the predator starvation rate in absence of prey \citep{baxter2008cost,sabo2005stochasticity}. The parameters $r, K, \alpha, \gamma$ and $c$ are real and positive constants. The system \eqref{eq:Baxter} has been used to calculate the annual population densities after a time step of one year \citep{baxter2008cost}. To support our investigations, we extend the Ricker\textendash type model \label{Baxter} with arbitrary step size, which is given by
\begin{gather}\label{DRM}
	\begin{split}
		N(t+h)&=N(t)\left( 1+h\left( e^{X(t)}-1\right) \right) \\
		P(t+h)&=P(t)\left( 1+h\left( e^{Y(t)}-1\right) \right), 
	\end{split}
\end{gather}
where
\begin{gather}\label{x_y}
	\begin{split}
		X(t)&=r(1-\frac{N(t)}{K})-\alpha P(t) \\
		Y(t)&=\alpha \gamma N(t)-c
	\end{split}
\end{gather}
and $h$ is a constant step size. The system \eqref{DRM} is the generalised version of \eqref{eq:Baxter} that considers the unit increments as a parameter such that the system \eqref{DRM} is the same as system \eqref{eq:Baxter} if $h=1$. 

In this paper, we also look into another popular population model, the model of Lotka\textendash Volterra\textendash. The generalised version of the discrete Lotka\textendash Volterra model is given by 
\begin{gather}\label{DVLM}
	\begin{split}
		N(t+h)&=N(t)\left(1+hX(t)\right)\\
		P(t+h)&=P(t)\left(1+hY(t)\right), 
	\end{split}
\end{gather}
where $h$ is the discrete step size.  In both forms of generalised discrete systems, the step size $h$ plays a critical role by permitting the users to choose appropriate time discretisation for each model. We observe that these generalized discrete systems  incorporating an arbitrary fixed step size are forward Euler's approximations of the respective continuous\textendash time systems.

Beyond the step size selections, the robustness of model parameters is critical, however, it is often a source of uncertainty in models based on real data. A slight variation of model parameters may change the equilibria and directly affect the system stability and robustness of solutions of the system \citep{enatsu2012global,jana2013chaotic,wang2020reinforcement}. If the parameters are estimated from data, lack of information and the inability to collect sufficient real\textendash world data in ecological systems can lead to an imprecise set of model parameters \citep{hines2014determination}. Therefore, investigating a suitable set of parameters that agrees with the selection of stable or unstable dynamics is essential when constructing population models, especially in approximating continuous\textendash time systems. 

In addition, discrete systems derived from first principles have common proprieties as the discretised approximations of the continuous\textendash time systems if some conditions are satisfied. The correspondence between the discrete\textendash time model and the continuous\textendash time model is the discrete mapping where step size is treated as a parameter. This builds a platform to link both discrete and continuous systems to analyse further the interaction between step size $h$ and the model parameters that affect the model performance and stability. 

Following the idea of parameterising the step size, \citep{liu2007complex} has derived the stability properties for the continuous\textendash time Lotka\textendash Volterra type predator\textendash prey system with scaled model parameters, and showed unstable and stable population dynamics for derived conditions in step size selections (see \citep{din2019stability,krivine2007discrete} for similar studies for Lotka\textendash Volterra type predator\textendash prey models). To the best of our knowledge, no study has investigated the dynamic inconsistency under arbitrary step size for the Ricker\textendash type ordinary differential equation(ODE) predator\textendash prey model \citep{brauer2012mathematical}, given by   
\begin{gather}\label{RODEM}
	\begin{split}
		N'(t)&=N(t)\left( e^{X(t)}-1\right)\\
		P'(t)&=P(t)\left( e^{Y(t)}-1\right).
	\end{split}
\end{gather} 
Hence, we study the required conditions for stable and unstable population dynamics of ODE system \eqref{RODEM} and their discretised system \eqref{DRM} with generalised step size. Consequently, the results identify similar or different dynamical properties of approximated discrete systems compared to the corresponding continuous\textendash time model. 

We perform a comparable study on qualitative analysis of the stability properties of a commonly used continuous\textendash time population model, a logistic growth Lotka\textendash Volterra type predator\textendash prey system \citep{zhao2020complexity}
\begin{gather}\label{VLODEM}
	\begin{split}
		N'(t)&=N(t)X(t)\\
		P'(t)&=P(t)Y(t).
	\end{split}
\end{gather}
Here, the prey population is influenced by prey natural growth, prey restricted growth in terms of prey carrying capacity, and prey death caused by predator attacks. The predator abundance is governed by population growth due to predation and natural death, respectively. We indicate that the generalised Lotka\textendash Volterra model \eqref{DVLM} is the approximated discrete system of the continuous\textendash time ODE model \eqref{VLODEM}. 

The dynamical properties of predator\textendash prey systems of continuous\textendash time models have been studied extensively in the literature (see examples at \citep{windarto2020modification,ackleh2015competitive,merdan2010stability}), but the factors that can perturb the dynamical properties (e.g. stabilising or destabilising the system) has not been fully analysed for discrete approximations of continuous\textendash time systems. We will fill the gap of identifying the constraints that deliver similar or different dynamical properties of two continuous\textendash time models and their approximated discrete system in terms of arbitrary step size under parameter space. More precisely, we derive the stability properties of discrete systems associated with corresponding continuous\textendash time models such that the interconnection of these two systems is purely observable. The comparison of both approximated (discrete) and actual (continuous) systems reveals the conditions that must be satisfied by the time discretisation. Our study also highlights the importance of the changes in (some) model parameters that direct the system to have stable or unstable population dynamics. Therefore, our work contributes to choosing the best\textendash suited step size for a particular discrete\textendash time system and helps to understand stabilising and destabilising factors of the continuous system and approximated discretised system under parameter space.  

This paper is organised as follows. In Section \ref{sec:chap3_stability}, a qualitative study to construct stability constraints is carried out through a Jacobian analysis of the discrete systems followed by that of the continuous systems. In Section \ref{sec:connection}, we derive the dynamical properties of discrete systems connected to the continuous\textendash time systems and investigate the conditions that lead the system solutions to become stable in terms of selecting suitable sets of model parameters and determining the effects of the step size in time discretisation. We then demonstrate population dynamics to justify our theoretical findings through numerical simulations in Section \ref{sec:chap3_results}.
\section{Equilibrium stability and dynamics}\label{sec:chap3_stability}

To understand the interconnection of dynamical properties in generic discretised systems with continuous\textendash time systems under arbitrary step size, we investigate the equilibrium points and the possible stable states of all Ricker\textendash type and Lotka\textendash Volterra systems mentioned in this paper. We first examine the system dynamics at the fixed points of the generalised discrete systems \eqref{DRM} and \eqref{DVLM}. 

In order to develop stability conditions for discrete systems, first consider the fixed point iteration formulation of a nonlinear map:
\begin{equation}\label{LC_Discrete}
	\chi(t+h)=F(\chi(t)),
\end{equation}
where $\chi(t)\in\mathbb{R}^m$ and $F$ has a Lipschitz condition. A point $\eta$ is said to be a fixed point satisfying $\eta=F(\eta)$ where $F$ is a map such that $F:I\rightarrow I$ and $I$ is a region in $\mathbb{R}^m$. It is proven that the fixed point $\eta$ is asymptotically stable if there exists a norm such that
\begin{equation}\label{discretestabilitycondition}
	||J_F||<1,
\end{equation} 
where $J_F$ is the Jacobian matrix of the discrete system evaluated at $\eta$ \citep{alligood1996two}. Since we consider two dimensional systems we have the following characterisations at $\eta$ in terms of the two eigenvalues $\lambda_1$ and $\lambda_2$ of the Jacobian matrix of the iterated map:

\begin{enumerate}
	\item[(i)] $\left|\lambda_i \right|<1,\>i=1, 2 $; a sink, locally asymptotically stable
	\item[(ii)] $\left|\lambda_i \right|>1,\> i=1, 2 $; a source
	\item[(iii)] one of $\left|\lambda_i \right|>1$ and other $\left|\lambda_i \right|<1,\> i=1, 2 $; a saddle
	\item[(iv)] one of $\left|\lambda_i \right|=1$ and  other $\left|\lambda_i \right|\neq1,\> i=1, 2$; non\textendash hyperbolic.
\end{enumerate}
Note that our investigation is based on this eigenvalue classification. Also, we can get flip, saddle and Hopf bifurcations when we traverse the boundaries of different stability domains. 

The equilibrium points for the Ricker\textendash type and Lotka\textendash Volterra discrete models can be obtained by solving 
\begin{gather}\label{discrete_N_P}
	\begin{split}
		N(t+h)&=N(t)\\
		P(t+h)&=P(t)
	\end{split}
\end{gather}
for $N$ and $P$. For the generalised Ricker\textendash type discrete model, \eqref{discrete_N_P} implies $X(t)=Y(t)=0$ and this also holds for the discrete Lotka\textendash Volterra model, \eqref{DVLM}. The fixed points of both the discrete formulations of the generalised Ricker\textendash type discrete model and the discrete Lotka\textendash Volterra model (\eqref{DRM} and \eqref{DVLM}, respectively) are the same, namely, $E_1\equiv(0,0)$, $E_2\equiv(K, 0)$ and $E_3\equiv\left (\frac{c}{\alpha \gamma}, \frac{r}{\alpha}\left ( 1-\frac{c}{K\alpha \gamma} \right )\right )$. Further, these are the same equilibrium points obtained for the continuous time models since the same condition $X(t)=Y(t)=0$ was satisfied, and we use the same notation. We denote the Jacobians of discrete maps in \eqref{DRM} and \eqref{DVLM} by $\widehat{J_R}$ and $\widehat{J_{LV}}$, respectively. Then some analysis gives, after omitting the dependence on $t$,
\begin{equation*}
	\widehat{J_R}=\begin{bmatrix}
		1+h\left( e^X-1-\frac{r}{K}Ne^X\right)  & -\alpha Nhe^X\\ \alpha \gamma Phe^Y& 1+h\left( e^Y-1\right) 
	\end{bmatrix}
\end{equation*}
and
\begin{equation*}
	\widehat{J_{LV}}=\begin{bmatrix}
		1+h\left( X-\frac{r}{K}N \right) & -\alpha hN\\ \alpha \gamma hP&1+hY
	\end{bmatrix}.
\end{equation*}
By using Jacobians, we derive the connection of the stability properties of the discrete systems to their respective ODE systems in the next section. 

We then classify the equilibria of the continuous\textendash time models \eqref{RODEM} (Ricker\textendash type) and \eqref{VLODEM} (Lotka\textendash Volterra) based on the stability since unstable and stable equilibria behave differently in population dynamics. Both systems \eqref{RODEM} and \eqref{VLODEM} have fixed points, namely $E_1\equiv(0,0)$, $E_2\equiv(K, 0)$ and $E_3\equiv\left (\frac{c}{\alpha \gamma}, \frac{r}{\alpha}\left ( 1-\frac{c}{K\alpha \gamma} \right )\right )$. Then the Jacobian matrix of \eqref{RODEM} is
\begin{equation*}
	J_R=\begin{bmatrix}
		e^X-1-\frac{r}{K}Ne^X & -\alpha Ne^X\\ \alpha \gamma Pe^Y& e^Y-1
	\end{bmatrix}
\end{equation*}
and the Jacobian matrix for \eqref{VLODEM} is 
\begin{equation*}
	J_{LV}=\begin{bmatrix}
		X-\frac{r}{K}N & -\alpha N\\ \alpha \gamma P&Y
	\end{bmatrix}.
\end{equation*}	
We observe similar stability properties for both ODE models even though the Jacobian matrices are different. The Jacobian matrices evaluated at $E_1=(0,0)$ are	
\[J_R\mid_{E_{1}}=\begin{bmatrix}
	e^r-1 & 0\\ 
	0 & e^{-c}-1
\end{bmatrix} \text{and } 
J_{LV}\mid_{E_{1}}=\begin{bmatrix}
	r & 0\\ 
	0 & -c
\end{bmatrix}.\]
Thus the eigenvalues at $E_1$ are 
\begin{gather}\label{lambdaE1ODE}
	\begin{split}
		\lambda_R\mid_{E_{1}}&=\left\lbrace e^r-1, e^{-c}-1\right\rbrace\\
		\lambda_{LV}\mid_{E_{1}}&=\left\lbrace r,-c \right\rbrace
	\end{split}
\end{gather}
where $\lambda_R$ and $\lambda_{LV}$ are eigenvalues for the system \eqref{RODEM} and \eqref{VLODEM}, respectively. $E_1$ is unstable for both models regardless of any choices of parameter values since it is a saddle point such that one eigenvalue is positive and one is negative. This means that the population never returns to $E_1$ after a small deviation of the population variation. 

The Jacobian matrices evaluated at $E_2=(K,0)$ are	  
\[J_R\mid_{E_{2}}=\begin{bmatrix}
	-r & -K\alpha\\ 
	0 & e^{\alpha\gamma K-c}-1
\end{bmatrix}\text{and }
J_{LV}\mid_{E_{2}}=\begin{bmatrix}
	-r & -K\alpha\\ 
	0 & \alpha\gamma K-c
\end{bmatrix}\]
and so 
\begin{gather}\label{lambdaE2ODE}
	\begin{split}
		\lambda_R\mid_{E_{2}}&=\left\lbrace -r, e^{\theta}-1\right\rbrace\\
		\lambda_{LV}\mid_{E_{2}}&=\left\lbrace -r,\theta \right\rbrace
	\end{split}
\end{gather}
where $\theta=\alpha\gamma K-c$. The $\theta$ value indicates a condition to determine the properties of eigenvalues. Therefore, investigations for stability properties are presented in terms of $\theta$ where appropriate. Then, $E_2$ is asymptotically stable for both models if the parameters satisfy $\theta<0$ (i.e. $\alpha\gamma K<c$) since both eigenvalues are then real and negative. This means that if there is a small deviation of the population densities away from $E_2$, both prey and predators can return to prey\textendash carrying capacity ($N=K$) and no predators, respectively. This implies that the prey population can be sustained by reaching its maximum carrying capacity without having predators in the ecosystem even if a small number of predators return, or prey experiences high mortality. Moreover, this case shows the predator\textendash prey existence under low predator populations where their ability to survive through food availability is determined by $K<\frac{c}{\alpha\gamma}$. On the other hand, for both models, $E_2$ is an unstable saddle point if $\theta>0$, and if $\theta=0$ $E_2$ is a non\textendash hyperbolic point where the system dynamics depends on the nonlinear terms of the model equations since it cannot be predicted from the eigenvalue analysis of the Jacobian matrix. 

The Jacobian matrices evaluated at $E_3=\left (\frac{c}{\alpha \gamma}, \frac{r}{\alpha}\left ( 1\frac{c}{K\alpha \gamma} \right )\right )$ are 
\[J_R\mid_{E_{3}}=J_{LV}\mid_{E_{3}}=\begin{bmatrix}
	-\frac{rc}{\alpha\gamma K} & -\frac{c}{\gamma}\\ 
	r\gamma\left ( 1-\frac{c}{\alpha\gamma K} \right ) & 0
\end{bmatrix},\]
and since these Jacobian matrices are the same for both ODE models we obtain identical results for stability analysis. The eigenvalues calculated for $E_3$ satisfy the characteristic polynomial
\begin{equation}\label{stability_lambda}
	\lambda^2+T\lambda+D=0
\end{equation}
where $T=\frac{rc}{\alpha \gamma K}$ and $D=c\left(r-T\right)=\theta T$. Thus,
\begin{equation}\label{lambdaE3ODE}
	\lambda_R\mid_{E_{3}}=\lambda_{LV}\mid_{E_{3}}=\left\lbrace \lambda_1, \lambda_2 \right\rbrace
\end{equation}
where $\lambda_j=\frac{-T\pm\sqrt{T^2-4D}}{2}=\frac{-T\pm\sqrt{T\left( T-4\theta\right) }}{2}$, $j=\{1,2\}$. Note that if $T\in(0,4\theta)$, the imaginary component for the continuous\textendash time models is $\frac{\sqrt{T(T-4\theta)}}{2}i$ and the larger the imaginary component the more oscillatory are the dynamics (Note the maximum imaginary component occurs when $T=2\theta$). The system stability status at a particular equilibrium point can be observed by looking at the sign of the eigenvalues and whether they are real or complex. We can conclude that $E_3$ is asymptotically stable if 
\begin{equation}\label{stability_theta}
	\theta>0,
\end{equation} 
with oscillatory dynamics if $T\in (0,4\theta)$, and $E_3$ is an unstable saddle point if $\theta<0$. $E_3$ is the only non\textendash trivial equilibrium point that has positive populations for both prey and predators. Note that if $\theta=0$ we can have a non\textendash hyperbolic property.

\section{Deriving connections of dynamical properties in discrete and continuous systems}\label{sec:connection}

The connections between discrete systems and ODE systems are derived under arbitrary step size. Stability analysis for the discrete\textendash time systems is simplified using derivations from the respective ODE systems where necessary. These stability constraints are observed through a Jacobian analysis, and the factors that could affect system dynamics are analyzed through variations of model parameters and by defining a particular range of step size.

We observe that
\begin{equation*}
	\widehat{J_R}=I+hJ_R
\end{equation*}
and
\begin{equation*}
	\widehat{J_{LV}}=I+hJ_{LV}.
\end{equation*} 
Hence $\lambda(\widehat{J_R})=1+h\lambda(J_R)$ and $\lambda(\widehat{J_{LV}})=1+h\lambda(J_{LV})$, where the eigenvalues of $J_R$ and $J_{LV}$ are given in \eqref{lambdaE1ODE}, \eqref{lambdaE2ODE} and \eqref{lambdaE3ODE}. Thus from the eigenvalue analysis, we have locally asymptotic stability for discrete mappings if $|1+h\lambda(J_R)|<1$ and similarly for $\lambda(J_{LV})$. In the case that  $\lambda(J_{R}), \lambda(J_{LV})<0$ this lead to $-h\lambda(J_R)<2$ and $-h\lambda(J_{LV})<2$. Therefore, depending on the step size and the eigenvalues of the continuous\textendash time system, new conditions exist in discretised systems for stable population dynamics. 

From the eigenvalue classification, $E_1$ is not asymptotically stable but there is a saddle point if $h(1-e^{-c})<2$ or $0<hc<2$ for Ricker\textendash type and Lotka\textendash Volterra discrete maps, respectively. For $\theta<0$, $E_2$ is asymptotically stable for the Ricker\textendash type discrete model, if $h<\left\lbrace \frac{2}{r},\frac{2}{1-e^\theta}\right\rbrace $, and for the discrete Lotka\textendash Volterra model, if $h<\left\lbrace\frac{2}{r},\frac{-2}{\theta}\right\rbrace $. Furthermore, if $\theta=0$ and $h\neq\frac{2}{r}$, $E_2$ is non\textendash hyperbolic for both models. For $\theta<0$, $E_2$ is a saddle point for the Ricker\textendash type discrete model if one of the following conditions holds
\begin{enumerate}
	\item[(i)] $\theta < -r \textrm{ and } \frac{2}{r}< h<\frac{2}{1-e^\theta},$
	\item[(ii)] $\theta > -r \textrm{ and } \frac{2}{r}> h>\frac{2}{1-e^\theta}$. 
\end{enumerate}
For $\theta<0$, $E_2$ is a saddle point for the Lotka\textendash Volterra discrete model if one of the following conditions hold
\begin{enumerate}
	\item[(i)] $\theta < ln(1-r) \textrm{ and } \frac{2}{r}< h<\frac{-2}{\theta},$
	\item[(ii)] $\theta > ln(1-r) \textrm{ and } \frac{2}{r}> h>\frac{-2}{\theta}.$
\end{enumerate}
In the case of the fixed point $E_3$, we can classify the stability associated with the two cases if the eigenvalues of $J_R$ and $J_{LV}$ are given as $\lambda_1$ and $\lambda_2$. From \eqref{stability_lambda}, the eigenvalues satisfy $\lambda^2+T\lambda+D=0$ where $T=\frac{rc}{\alpha\gamma K}$ and $D=c(r-T)=\theta T$.  For $\theta>0$, 

\begin{enumerate}
	\item[(i)] if $\lambda_1$ and $\lambda_2$ are real and negative where $\lambda_1=\lambda_2$ then $E_3$ is asymptotically stable if \[h<\frac{T}{4}.\] This happens only if $T=4\theta$.
	\item[(ii)] if both $\lambda_1$ and $\lambda_2$ are real and negative where $\lambda_1\neq\lambda_2$ then $E_3$ is asymptotically stable if the step size satisfies 
	\[h<\left\lbrace-\frac{2}{\lambda_1},-\frac{2}{\lambda_2} \right\rbrace. \] 
	This happens if $T-4\theta>0$. 
	\item[(iii)] if both $\lambda_1$ and $\lambda_2$ are complex conjugate eigenvalues, say $a\pm ib$, then from the above $a^2+b^2=D=\theta T$ and $T=-2a$. This occurs when $T-4\theta<0$. With these complex eigenvalues, the population dynamics lead to oscillations with time. Then the bound for the step size is
	\begin{equation}\label{Discrete_condition}
		h<\left\lbrace\frac{-2a}{a^2+b^2} \right\rbrace=\frac{T}{D}=\frac{1}{\theta}.
	\end{equation}
	This can only happen if $0<1+hT(h\theta-1)$. Note that if $h=\frac{1}{\theta}$ then both eigenvalues have magnitude one.	    
\end{enumerate}

For $\theta<0$, say $\lambda_1>0$ and $\lambda_2<0$, then $E_3$ is non\textendash hyperbolic if $h=-\frac{2}{\lambda_2}$. For $\theta=0$, $E_3$ is non\textendash hyperbolic if $h\neq\frac{2}{T}$.

\begin{table}[H]
	\caption{Stability status for discrete and continuous Ricker\textendash type (RK) and Lotka\textendash Volterra (LV) models where $E_1\equiv(0,0)$, $E_2\equiv(K, 0), E_3\equiv\left (\frac{c}{\alpha \gamma}, \frac{r}{\alpha}\left ( 1-\frac{c}{K\alpha \gamma} \right )\right )$, and $\lambda_1, \lambda_2$ are eigenvalues of $E_3$ calculated from \eqref{lambdaE3ODE}.}
		\newcolumntype{C}{>{\centering\arraybackslash}X}
		\begin{tabular}{@{}cccccc@{}}
			\toprule
			\multirow{2}{*}{\textbf{Stability}}      & \multirow{2}{*}{\textbf{Model}} & \multicolumn{2}{c}{\textbf{$E_1$}}       & \multicolumn{2}{c}{\textbf{$E_2$}}  \\ \cmidrule{3-6}  
			&  & \textbf{Discrete} & \textbf{Continuous} & \textbf{Discrete} & \textbf{Continuous} \\ \cmidrule{1-6}
			\multirow{2}{*}{\textbf{\begin{tabular}[c]{@{}c@{}}Asym.\\ stable\end{tabular}}} & \textbf{RK} & - & - &if $\theta<0$ and $h<\{\frac{2}{r},\frac{2}{1-e^\theta}\}$& if $\theta<0$\\\cmidrule{2-6} & \textbf{LV}& -& -& if $\theta<0$ and $h<\{\frac{2}{r},-\frac{2}{\theta}\}$& as above\\ \midrule
			\multirow{2}{*}{\textbf{\begin{tabular}[c]{@{}c@{}}Non\textendash\\ hyperbolic\end{tabular}}} & \textbf{RK}& if $h=\frac{2}{1-e^{-c}}$ & -& \begin{tabular}[c]{@{}c@{}}if $\theta=0, h\neq\frac{2}{r}$ \\ if $\theta<0, h=\frac{2}{1-e^\theta}, h\neq\frac{2}{r}$ \\ if $h=\frac{2}{r}, \theta\neq0, \theta\neq ln(1-r)$\end{tabular}& if $\theta=0$\\\cmidrule{2-6}
			& \textbf{LV}& if $h=\frac{2}{c}$& -& \begin{tabular}[c]{@{}c@{}}if $\theta=0, h\neq\frac{2}{r}$ \\ if $\theta<0, h=-\frac{2}{\theta}, \theta\neq 0, \theta\neq -r$ \\ if $h=\frac{2}{r}, \theta\neq 0, \theta\neq -r$\end{tabular}& as above\\   \midrule
			\multirow{2}{*}{\textbf{Saddle}} & \textbf{RK}& if $h<\frac{2}{1-e^{-c}}$ & \begin{tabular}[c]{@{}c@{}}always a \\ saddle point\end{tabular}& \begin{tabular}[c]{@{}c@{}}if $\theta<0, \theta < -r, \frac{2}{r}< h<\frac{2}{1-e^\theta}$\\if $\theta<0, \theta > -r, \frac{2}{r}> h>\frac{2}{1-e^\theta}$\\ if $\theta>0, h<\{\frac{2}{r},\frac{2}{1-e^\theta}\}$\end{tabular}& if $\theta>0$ \\\cmidrule{2-6}
			& \textbf{LV}&if $h<\frac{2}{c}$ & \begin{tabular}[c]{@{}c@{}}always \\ a saddle point\end{tabular}& \begin{tabular}[c]{@{}c@{}} if $\theta<0, \theta < ln(1-r), \frac{2}{r}< h<-\frac{2}{\theta}$\\ if $\theta<0, \theta > ln(1-r), \frac{2}{r}> h>-\frac{2}{\theta}$\end{tabular}& as above\\
			\bottomrule    		                 
		\end{tabular}
		\begin{tabular}{@{}cccc@{}}
			\toprule
			\multirow{2}{*}{\textbf{Stability}} & \multirow{2}{*}{\textbf{Model}} & \multicolumn{2}{c}{\textbf{$E_3$}} \\ \cmidrule{3-4}  
			& & \textbf{Discrete} & \textbf{Continuous}  \\ \cmidrule{1-4}
			\multirow{2}{*}{\textbf{\begin{tabular}[c]{@{}c@{}}Asym.\\ stable\end{tabular}}} & \textbf{RK} & \begin{tabular}[c]{@{}c@{}}if $\theta>0$, $T=4\theta, 0<h<\frac{T}{4}$\\if $\theta>0, T-4\theta>0, h<\{-\frac{2}{\lambda_1}, -\frac{2}{\lambda_2}\}$\\if $\theta>0, T-4\theta<0, h<\frac{1}{\theta}, 0<1+hT(h\theta-1)$\end{tabular}&  if $\theta>0$ \\\cmidrule{2-4} & \textbf{LV}& as above& as above\\ \midrule
			\multirow{2}{*}{\textbf{\begin{tabular}[c]{@{}c@{}}Non\textendash\\ hyperbolic\end{tabular}}} & \textbf{RK}& \begin{tabular}[c]{@{}c@{}}if $\theta>0, T-4\theta>0, h=-\frac{2}{\lambda_i}, h\neq-\frac{2}{\lambda_j}, i\neq j, i,j=\{1,2\}$ \\ if  $\theta<0, h=-\frac{2}{\lambda_2}, \lambda_2<0, \lambda_1>0$ \\ if $\theta=0, h\neq \frac{2}{T}$\end{tabular}& if $\theta=0$\\\cmidrule{2-4}
			& \textbf{LV}& as above & as above\\   \midrule
			\multirow{2}{*}{\textbf{Saddle}} & \textbf{RK} & if $\theta>0, T-4\theta>0, -\frac{2}{\lambda_i}<h<-\frac{2}{\lambda_j}, i\neq j, i,j=\{1,2\}$ & if $\theta<0$ \\\cmidrule{2-4}
			& \textbf{LV}& as above& as above\\\bottomrule    		                 
		\end{tabular}
	\label{tab:one}
\end{table}

Overall, the stability constraints are different in continuous\textendash time models and their corresponding discrete systems. A summary of the stability analysis of all eigenvalues for the Ricker\textendash type and Lotka\textendash Volterra discrete and continuous\textendash time models are given in Table \ref{tab:one}. Stability criteria evaluated at equilibrium point $E_3$ are similar for both Lotka\textendash Volterra and Ricker\textendash type models. At equilibrium point $E_2$, the stability conditions are similar for both models under continuous\textendash time setting only, and different otherwise. We consider the case (iii) when $\theta>0$ for further simulations since it is a stable spiral where the population returns to a steady state, $E_3$. Thus, additional constraints are required for stable population dynamics in approximated discrete systems than the continuous\textendash time system (which is $\theta>0$). Therefore, the model dynamics of approximated discrete systems depends on the selected step size and (some) model parameters. We will use \eqref{Discrete_condition} to observe the stability with different step sizes and ranges of model parameters in the next section.

\section{Numerical results}\label{sec:chap3_results}

We present a numerical simulation study for the Lotka\textendash Volterra and Ricker\textendash type discrete systems to illustrate the theoretical findings discussed in Section \ref{sec:connection}. The stability condition \eqref{Discrete_condition}, where the discrete systems become stable at $E_3$, are demonstrated for some parameter ranges and step sizes. The impact of model parameter variations on the system that changes the model dynamics are then investigated through a few examples. This numerical study clearly demonstrates the changes to the population over time in the long\textendash term scale.

\begin{figure}[H]
	\centering		
	\includegraphics[width=13cm]{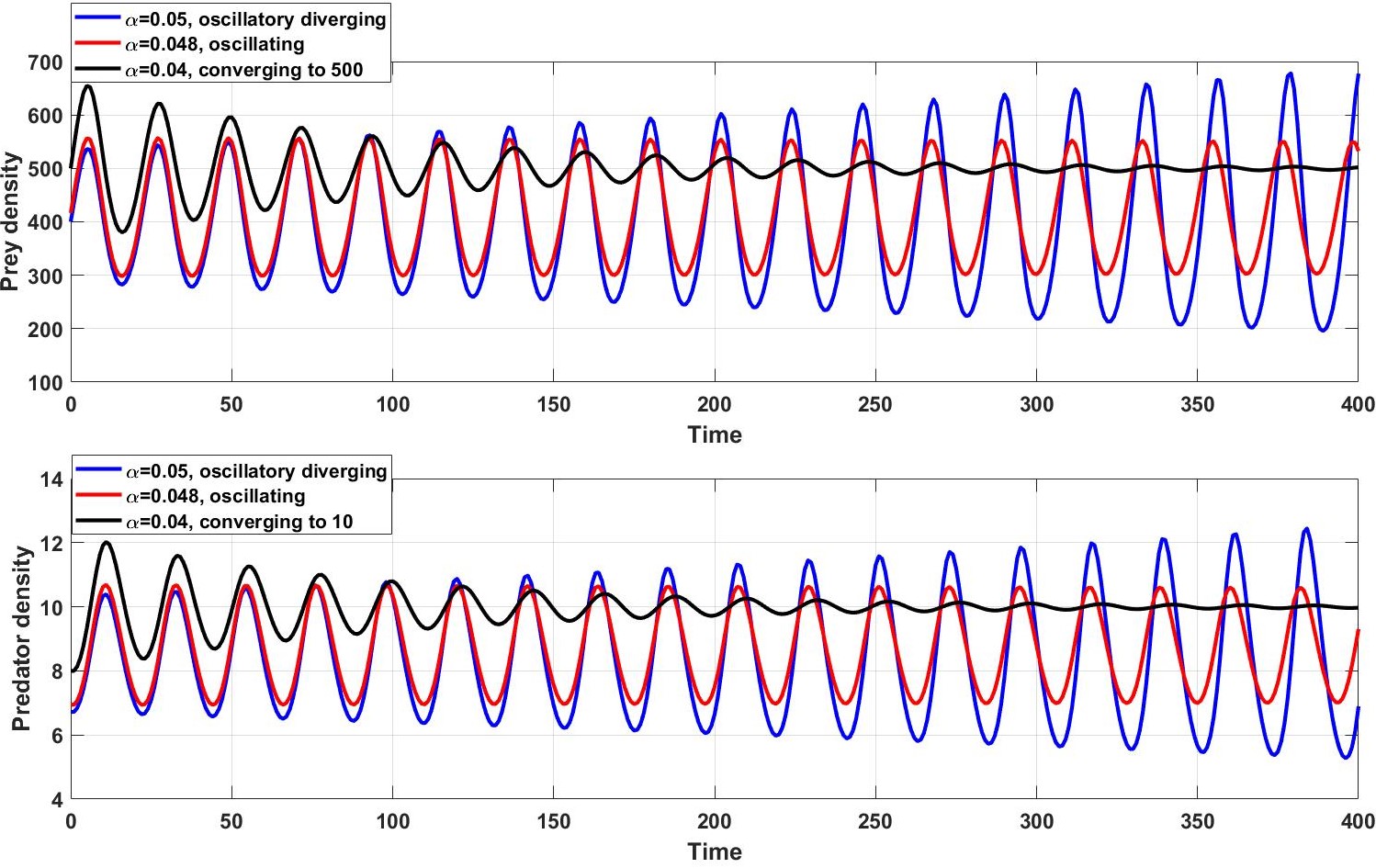}		
	\caption{Three different dynamics of discrete Ricker\textendash type system \eqref{Baxter} with slightly varying $\alpha$ where $K=2500, \gamma=0.01, c=0.2$ and $r=0.5$. For different $\alpha=\{0.05, 0.048, 0.04\}$ values, predator\textendash prey populations are diverging, seeming to converge, and converging, respectively. Here, the system \eqref{Baxter} is derived for a unit step size, which is similar to the system \eqref{DRM} when $h=1$.}
	\label{fig:1_dynamic}
\end{figure}

We observed that some parameters impact the system dynamics and can stabilise or destabilise the populations if the time discretisation is fixed. In our discrete\textendash time models, model stability is governed by the parameters $\alpha, \gamma, K, c$ and $r$ according to the derivation of \eqref{Discrete_condition} given that the step size $h$ is fixed. We first observe the stability condition \eqref{Discrete_condition} for discrete systems at unit step size, $h=1$ by assigning the parameter values described in \citep{baxter2008cost}, where $\alpha=0.05, \gamma=0.01, K=2500, c=0.2, r=0.5$, then, $\theta=\alpha\gamma K-c=\frac{21}{20}, T=\frac{rc}{\alpha \gamma K}=\frac{2}{25}$. Then \eqref{Discrete_condition} specifies stability for step size $$h<\frac{20}{21},$$ where $\theta>0, T-4\theta<0$ and $1+hT(h\theta-1)>0.$ In this case, since the eigenvalues are complex, the size of the imaginary component of the eigenvalue evaluated at equilibrium point $E_3$ is $\frac{\sqrt{206}}{50}$. For this choice of parameters with step size $h=1$, the system is not asymptotically stable at $E_3$ since annual discretisation does not satisfy $h<\frac{20}{21}$. Thus, the system dynamics diverge while oscillating around the equilibrium point $E_3$ (blue curves in Figure \ref{fig:1_dynamic}). If now, $\alpha=0.048$ with the other parameters the same and a one year step size $h=1$, then the system may not become asymptotically stable since $\theta=1, T-4\theta<0$ and $1+hT(h\theta-1)>0$. In this case, the condition for stability in equation \eqref{Discrete_condition} is violated, $h=\frac{1}{\theta}=1$. The populations seem to be converging to $E_3$, but oscillate around $E_3$, (red curves in Figure \ref{fig:1_dynamic}). If $\alpha=0.04$ and $h=1$ with the other parameters the same, then $\theta=\alpha\gamma K-c=\frac{1}{20}>0, T-4\theta<0$ and $1+hT(h\theta-1)>0$, thus, $h<\frac{1}{\theta}$ and the condition \eqref{Discrete_condition} is true for this case. Therefore, the populations converge to $E_3=(500,10)$, see the black curves in Figure \ref{fig:1_dynamic}. Thus, there is considerable sensitivity to the choice of parameters and hence to the step size which controls the dynamics of the model.

\begin{figure}[H]
	\centering
	\includegraphics[width=13cm]{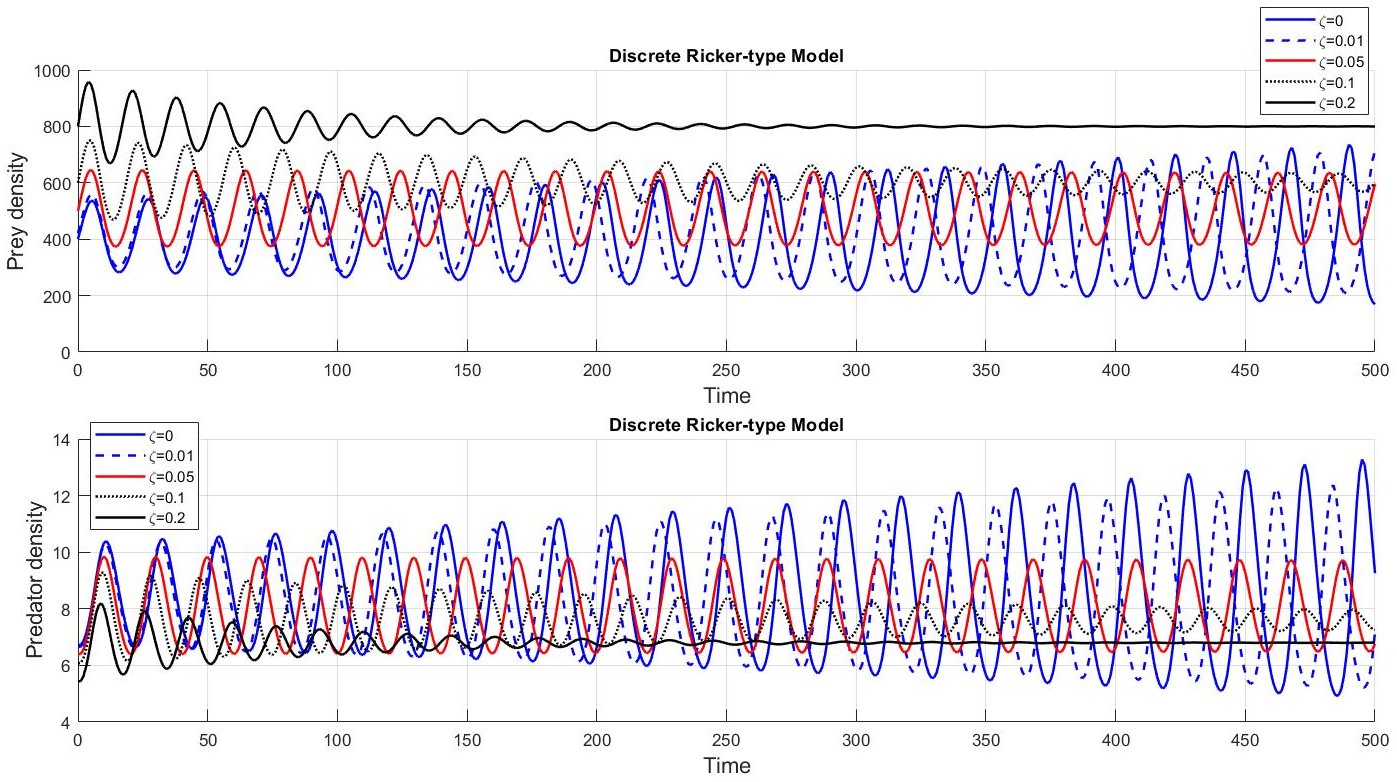}
	\includegraphics[width=13cm]{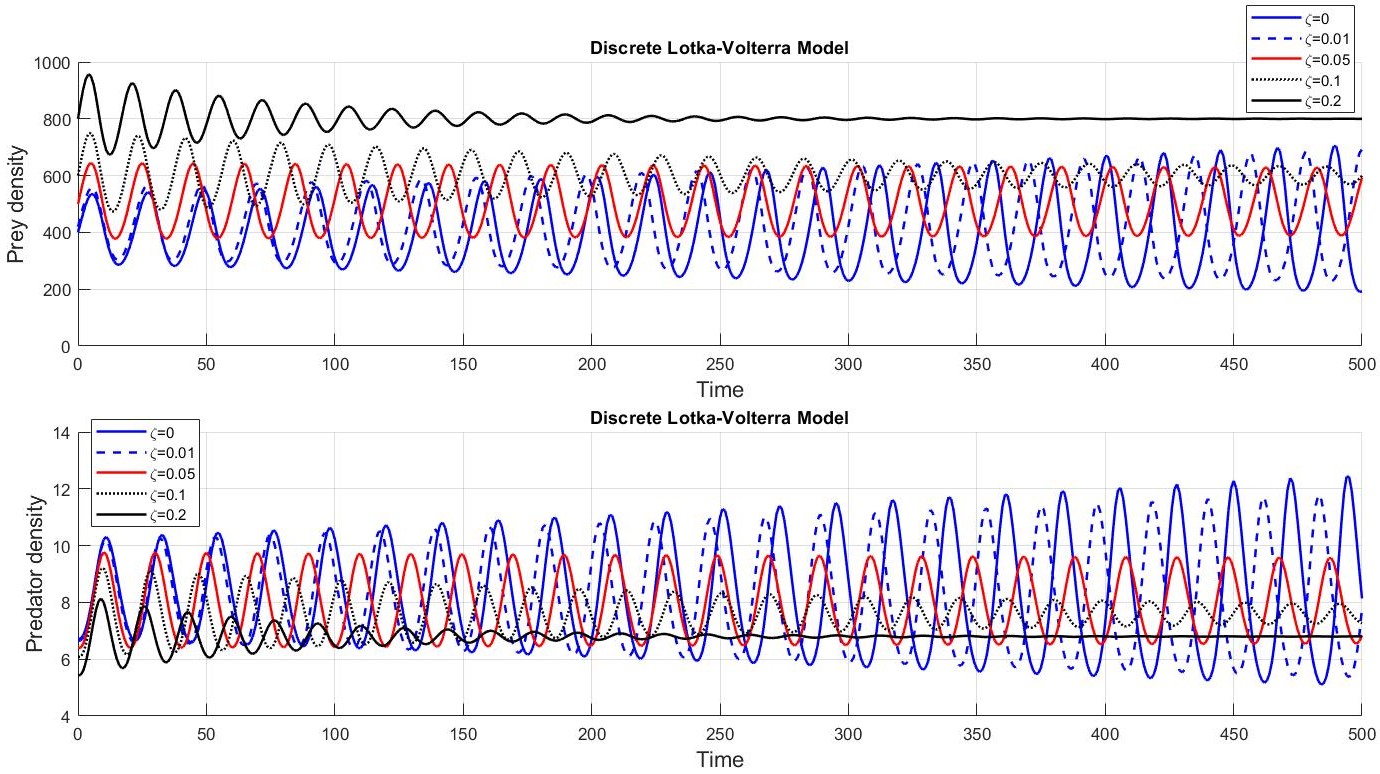}
	\caption{Different predator\textendash prey dynamics of discrete Ricker\textendash type and Lotka\textendash Volterra models with slightly varying parameter $c$ as $c+\zeta$ by $\zeta=\{0,0.01,0.05,0.1,0.2\}$ values where $h=1, K=2500, \gamma=0.01, \alpha=0.05, c=0.2$ and $r=0.5$. } \label{fig:4_dynamic}
\end{figure}

Moreover, if we choose $c=0.25+\epsilon$, rather than $c=0.2$ as previously, where $\epsilon$ is a small positive value and the other parameters are as in \citep{baxter2008cost} then 
$$\theta=\alpha\gamma K-c=\frac{25}{20}-\frac{5}{20}-\epsilon=1-\epsilon$$
and so \eqref{Discrete_condition} gives 
$$h<\frac{1}{1-\epsilon}.$$ 
Then for the step size $h=1$, the condition \eqref{Discrete_condition} is true such that the system generates populations that converge to the fixed points (except for $\epsilon=0$). Therefore, if the step size is fixed, a small deviation of a model parameter or specific parameters that define model stability can make dramatic changes in population dynamics.

Returning to the case of $c=0.2$, we varied $c$ by $c+\zeta$ for the step size $h=1$ where $\zeta$ is a small positive constant value. The numerical results of population densities over time in Figure \ref{fig:4_dynamic} confirm the stability conditions developed for discrete models as in equation \eqref{Discrete_condition}. For small values of $\zeta$ near zero, populations are oscillatory diverging, see the blue curves in Figure \ref{fig:4_dynamic}. A special dynamics is observed for $\zeta=0.05$ (red curves in Figure \ref{fig:4_dynamic}), in which case $c=0.25$ and corresponds to the case of $\zeta=0$ studied previously. Fixed point convergence is observed for $\zeta=\{0.1,0.2\}$ which satisfy the stability condition in \eqref{Discrete_condition}, see the black curves in Figure \ref{fig:4_dynamic}. The population dynamics according to the small changes of the parameter $c$ shows how critical the sensitivity of the parameters is. This is valid for other parameters that affect the stability of the discrete model, in our case $\alpha, \gamma, K, c$ and $r$. Moreover, the numerical simulations of the discrete models support the developed analytical results in discrete systems.  

\begin{figure}[H]
	\centering
	\includegraphics[width=8.5cm]{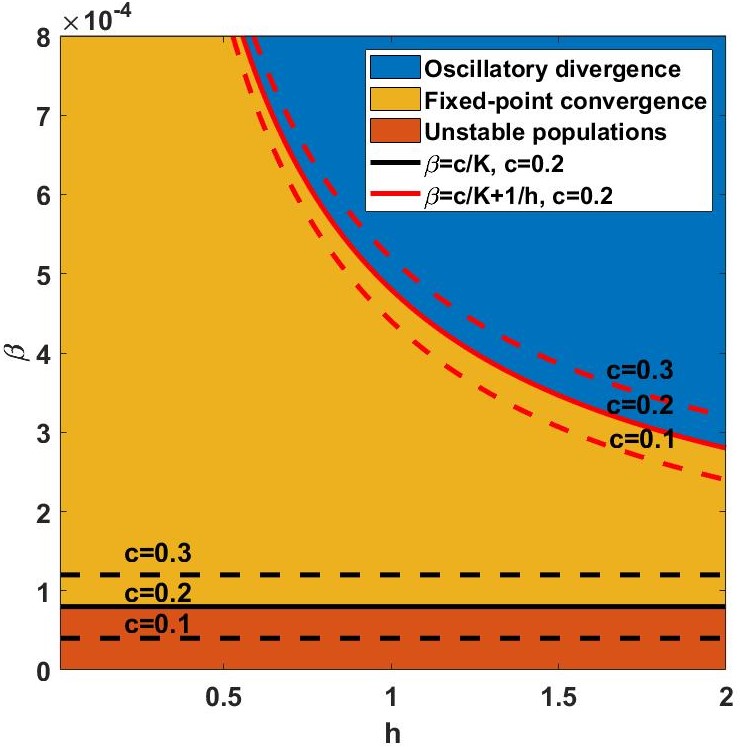}
	\caption{Stability regions of discrete Ricker\textendash type model and discrete Lotka\textendash Volterra model where $h$ is the step size, $\beta=\alpha\gamma$, $K=2500$ and $r=0.5$. The fixed\textendash point convergence region is bounded by $\beta=\frac{c}{K}+\frac{1}{h}$ and $\beta=\frac{c}{K}$, and boundary changes are marked in red and black lines for different $c=\{0.1, 0.2, 0.3\}$.  Here, $h$ is an independent variable and $\beta$ is calculated accordingly. If $\beta$ is known, the step size can be chosen to achieve the required stability status. Note that the stability regions are coloured for $c=0.2$, as represented in solid lines, and boundary shifting movements are displayed for $c=0.1$ and $c=0.3$ as represented in dashed lines. }
	\label{fig:2_dynamic}
\end{figure}

Figure \ref{fig:2_dynamic} shows the fixed\textendash point convergent region (yellow) for $\beta\in[0,8\times10^{-4}]$ where $\beta=\alpha\gamma$. Therefore, parameter $c$ has different boundaries for converging populations and preserving positive predator populations as demonstrated in the red and black curves of Figure \ref{fig:2_dynamic}, respectively. 

Special dynamics of the system \eqref{Baxter} is observed when $\alpha=0.048$ that is on the upper bound of the fixed\textendash point converging region (red curves in Figure \ref{fig:1_dynamic}). The black curves in Figure \ref{fig:1_dynamic} represent the converging predator and prey densities to their fixed points as the $\alpha$ is chosen from the yellow region of Figure \ref{fig:2_dynamic}. If $\alpha$ is chosen from the oscillatory divergence region, predator\textendash prey densities tend to oscillate continuously and increase in amplitude, see the blue curves in Figure \ref{fig:1_dynamic}.

\begin{figure}[H]
	\centering
	\includegraphics[width=8.5cm]{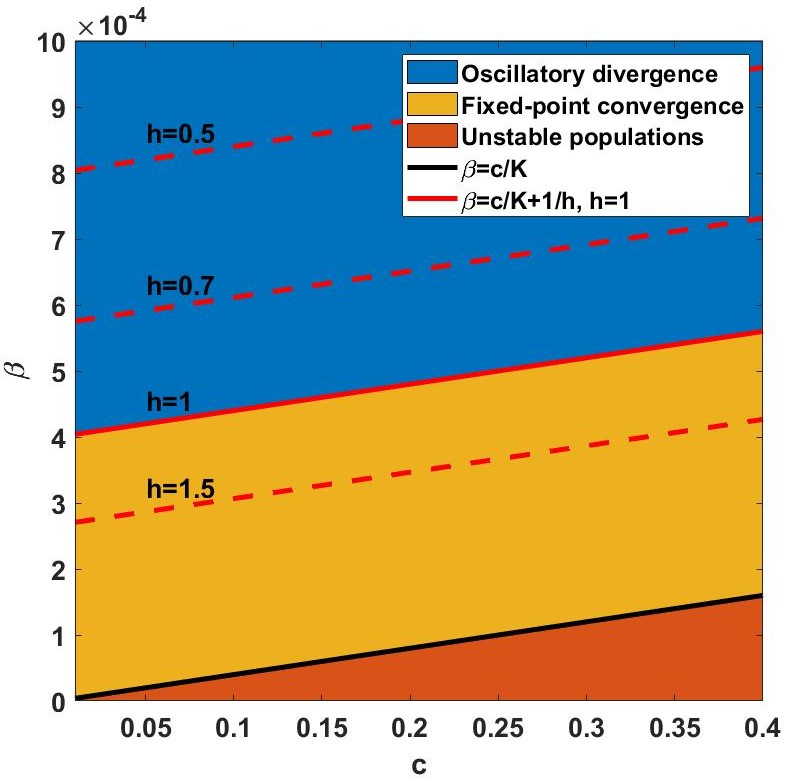}
	\caption{Expansion of fixed point convergence region with decreased step size ($h$) for discrete Ricker\textendash type model and discrete Lotka\textendash Volterra model where $\beta=\alpha\gamma$, $K=2500$ and $r=0.5$. The fixed\textendash point convergence region is bounded by $\beta=\frac{c}{K}+\frac{1}{h}$ and $\beta=\frac{c}{K}, \forall h>0$. Here, $c$ is the independent variable and $\beta$ is calculated accordingly. For $h=1$, the stability regions are coloured, and upper and lower boundaries of fixed\textendash point region are plotted for $\beta=\frac{c}{K}+1$ and $\beta=\frac{c}{K}$, as displayed in solid red and black lines, respectively. The upper boundary of fixed\textendash point region moves upward with decreasing step size, $h=1.5, h=0.7$ and $h=0.5$, and these upper boundary shifting movements are marked as red\textendash dashed lines. Note that the lower boundary of fixed\textendash point convergence region is valid for any $h$.}
	\label{fig:3_dynamic}
\end{figure}

On the other hand, from the theoretical analysis, we observed that the dynamical properties of the approximated discrete systems depends on a suitable time discretisation if the model parameters are fixed. As a guide to select a suitable step size $h$ with the required stability property, the impacts of variable step sizes are then investigated. Thus, different upper bounds on $h$ exist for the fixed point convergent regions (see Figure \ref{fig:3_dynamic}). The fixed\textendash point convergence region becomes larger for small step sizes. Figure \ref{fig:3_dynamic} indicates that larger step sizes are more likely to show oscillatory divergence. Populations become more stable and converge to the equilibrium point if the step size is small. This generates a time discretisation (into small time intervals) where the discrete systems behave more closely to the continuous systems. Certainly, to preserve the characteristics in numerical simulations, this idea supports the small step size recommendation in Euler's scheme, since it is a first\textendash order method \citep{efimov2019discretization}. These results lead to selecting suitable values for the step size in terms of preserving the required stability states when approximating discrete systems to respective continuous systems. On the other hand, if the step size is too small then on larger systems, the simulations may take a very long time to run, especially over large time intervals.  

\begin{figure}[H]
	\centering
	\includegraphics[width=13cm]{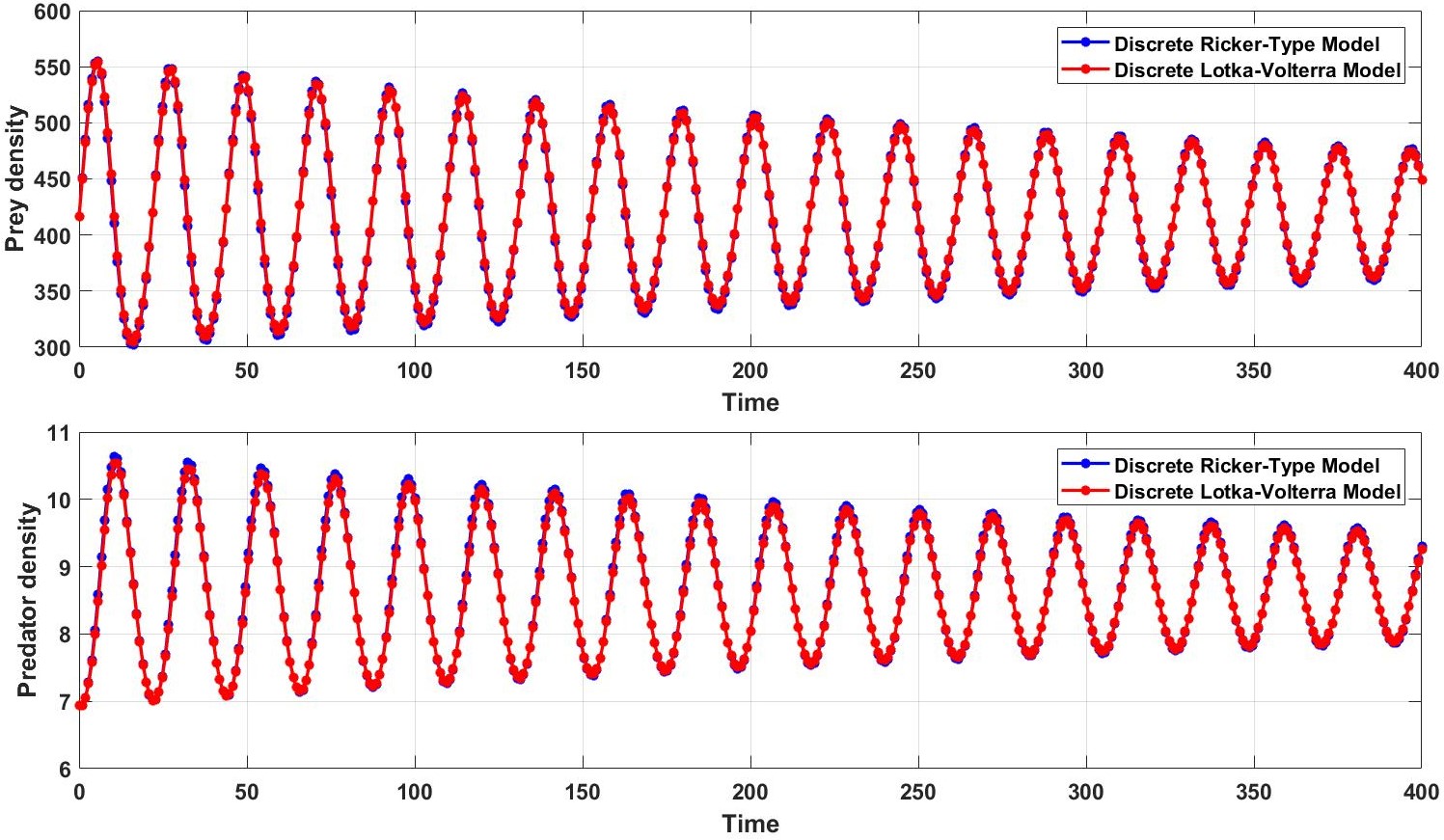}
	\caption{Special behaviour of predator\textendash prey populations for Ricker\textendash type and Lotka\textendash Volterra discrete models if $h=\frac{1}{\theta}=\frac{20}{21}$, where $K=2500, \gamma=0.01, \alpha=0.05, c=0.2$ and $r=0.5$. Predator\textendash prey populations seems to converge to a fixed point at the beginning, however, after a long\textendash time, the populations oscillate around the fixed point. Note that this exceptional case occurs only at the upper bound of the fixed\textendash point convergence region.}
	\label{fig:5_dynamic}
\end{figure}

Even though the model structures of the discrete Ricker\textendash type and Lotka\textendash Volterra models are different, the derived stability conditions for these models are similar in some cases. We examine the frequency plots of the population dynamics for both models with different step sizes. The red and blue curves in Figure \ref{fig:5_dynamic} provide evidence of the special behaviour when $\theta=\alpha\gamma K-c=21/20$ (with $h=\frac{1}{\theta}$). These curves have a similar pattern that converges to a fixed point at a later time. The continuously converging pattern of these curves show similar characteristics of the fixed\textendash point convergence at the beginning, and from the long\textendash time observations, the populations never approach fixed points. Therefore, this is an exceptional case that behaves as an upper bound for the fixed\textendash point convergence region. 

In summary, the stability criteria for the continuous\textendash time models depend only on the model parameters defined by $\theta$, however, when the continuous\textendash time models are discretised through numerical approximations, the stability criteria depend on the step size along with (some) model parameters. Therefore, the dynamical properties of the original ODE systems are different from the respective discretised systems with the same parameter values unless a suitable step size is defined. To get a better understanding of the model dynamics, the discretised systems should be thoroughly investigated under an eigenvalue analysis to identify suitable step sizes that agree with model parameters. In real predator\textendash prey systems, a prior knowledge of the system behaviour enhances the understanding of the factors that stabilise or destabilise the populations over time. If the parameters are estimated from a data set, this type of theoretical study provides a detailed analysis for selecting a suitable numerical simulation method that discretises the time component with a step size. 

\section{Discussion}

We have investigated the stability of the population dynamics of a Ricker\textendash type and Lotka\textendash Volterra discrete and continuous\textendash time population models. Based on the Jacobian analysis, important constraints are generated to identify the stability regions of the discrete systems. The generalised discrete systems can be viewed as the Euler numerical scheme of the approximate solution of the continuous\textendash time models. The discretised solution may or may not have the same properties as observed in the original continuous\textendash time model. Therefore, a qualitative analysis of the model dynamics in both the original continuous\textendash time model and the relevant discretised system is essential when modelling to inform ecological decisions. We showed that the two models have similar conditions to stabilise the system depending on the constraints as in Table \ref{tab:one}. This novel work increases the understanding of the similar behaviours of two structurally different predator\textendash prey models, in a nonlinear setting. 

The derived stability properties for continuous\textendash time and approximated discrete\textendash time solutions are different at each equilibrium point. Therefore, choosing a suitable time discretisation depends on the stability criteria that is determined by the dynamical properties. We assume the populations in real ecosystems coexist at, or move to, a stable\textendash state which is $E_3$ with asymptotic stability. Therefore, selecting a suitable time discretisation is essential to approximate the respective continuous\textendash time system if the model parameters are given.        

We highlighted the significance of understanding the choices of model parameters that impact the stability of the system. Some parameters need to be verified carefully due to the strong effect on system dynamics such as the parameters determined by $\theta$. Small changes in these parameters lead to large deviations in population count when the step size of the time discretisation is fixed (e.g. one year). A prior analysis on the impact of selecting a suitable time discretisation that is relevant to selected parameter sets is essential for a better performance of the models. We also numerically showed that the small changes in parameter values and step sizes stabilise or destabilise the system and form different dynamics in population densities. 

As we observed, parameter values can change the stability conditions of the model and result in abrupt dynamical behaviour in well\textendash formulated ecological systems. The factors that lead to the system becoming unstable could be useful in deciding population management plans for unstable populations. For the model parameters evaluated from a data set, a qualitative study that determines the stable or unstable dynamical properties is essential to understand the population dynamics. A long\textendash term stability analysis is required as some systems need population density variations over a long time duration due to long\textendash term impacts on populations. The effects of individual choices of parameters or suitable sets of parameters to preserve the stability of the systems must agree with an appropriate time discretisation in discretised systems.

We recommend performing a qualitative analysis on the dynamics of the desired approximation (discretised system) to identify the conflicts between actual (continuous\textendash time model) and approximated systems. Not limited to comparing the dynamical properties of approximated and actual systems under fixed parameter space, choosing a suitable step size also depends on the accuracy of the desired numerical approximation tool, such as the order of accuracy in numerical schemes \citep{islam2015comparative}. For some cases, the established non\textendash standard numerical approximations (e.g.\citep{mickens2005dynamic}) are suitable to answer the dynamic inconsistencies in approximated systems. Some discrete approximations generated through non\textendash standard discretisation methods behave as the original ODE system and dynamic consistency is independent of the chosen step size \citep{leslie1958stochastic,seno2007discrete}.

Our results on the Ricker\textendash type discrete model are consistent with the stabilisation that has been found in \citep{sabo2005stochasticity}. However, the discrete model stability analysis that shows the significance of dynamical properties over the parameter space goes beyond previous studies on investigating the stabilising and destabilising factors defined in \citep{din2013dynamics,merdan2009stability}. Besides, our results increase the understanding of how to select a suitable step size for the given model parameters under discrete settings, and when implementing a continuous\textendash time model through numerical simulations. 

We restricted our simulations to investigate the conditions that affect the system stability of population models. We carried out a theoretical study on classifying stability conditions for the discrete and continuous time models under consideration. Our numerical results are limited to demonstrate the properties only at $E_3$. We focused only on the Euler numerical scheme to generate discretised systems from continuous\textendash time models, however, other time discretisation methods could be considered. We showed that the stability of the discretised model depends on the selected step size if the model parameters are fixed. There are non\textendash standard discretisation methods of nonlinear ODEs to preserve the dynamic consistency regardless of the selection of the step size \citep{seno2007discrete,mickens2005dynamic}. We note that this is not a complete study of the dynamics of discrete mapping. For example from \eqref{LC_Discrete}, we note $\chi(t+2h)=F(\chi(t+h))=F(F(\chi(t)))$. Hence, we can study the fixed points of $F\circ F$ and the ensuing dynamics would lead to period two dynamics, depending on the nature of the map $F$. The two models will have different dynamics in this regard. The theories on limit cycles and bifurcation analysis \citep{rana2020chaotic,yousef2019stability,luis2011stability} were not considered, but, this could be seen through our numerical simulations. Of course, this can be generalised to period dynamics of any integer order and potentially lead to chaotic dynamics (of iteration of the logistic map).

The demonstrated stability analysis can be applied to other forms of two species continuous\textendash time and discrete\textendash time population models in returning more complex dynamical systems such as controlling species, functional responses, time delay and the Allee effect. Our work can be extended to study the dynamics of three or more species systems \citep{luis2017local} and to understand the stabilising and destabilising factors before obtaining the model outcomes. 

Finally, this work has implications in uncertainty quantification. In this setting, populations of models (with the same structure but different parameter sets) are constructed based on these models satisfying a set of common outputs. If there is sensitivity of the dynamics to the parameters, as is the case here, then this can make the process of uncertainty quantification also sensitive.

\section*{Author contributions}
All authors contributed equally.

\section*{Funding}
K.J.H. acknowledges support from the Australian Research Council Fellowship DE200101791.

\section*{Abbreviations}{
	The following abbreviations are used in this manuscript:\\
	
	\noindent 
	\begin{tabular}{@{}ll}
		ODE &Ordinary Differential Equations
	\end{tabular}
}

\bibliography{Ref}

\end{document}